\documentclass{amsart}
\input{psfig.sty}
\usepackage{amssymb}

\newtheorem{theorem}{Theorem}[section]
\newtheorem{lemma}[theorem]{Lemma}

\theoremstyle{definition}

\theoremstyle{remark}

\numberwithin{equation}{section}

\newcommand{\abs}[1]{\lvert#1\rvert}

\def\Z{\mathbb{Z}}

\def\Hc{Hall's conjecture }

\def\round#1{\left\lfloor#1\right\rceil}

\begin{document}

\title{A new algorithm to search for small nonzero $\abs{x^3-y^2}$ values}

\author{I. Jim\'enez Calvo}
%\address{C.I.S.C, C/Serrano 144, 28006--Madrid. (Spain).} 
%\email{ismaeljc@fdi.ucm.es} 
\address{C/Virgen de las Vi\~nas 11, 28031--Madrid. (Spain).} 
\email{ismaeljc@fdi.ucm.es,ijcalvo@gmail.com}

\author{J. Herranz}
\address{Dep. de Matem\`atica Aplicada IV,
Universitat Polit\`ecnica de Catalunya,
c/Jordi Girona, 1-3, 08034-Barcelona (Spain).}
\email{jherranz@ma4.upc.es} 

\author{G. S\'aez}
\address{Dep. de Matem\`atica Aplicada IV,
Universitat Polit\`ecnica de Catalunya,
c/Jordi Girona, 1-3, 08034-Barcelona  (Spain).}
\email{german@ma4.upc.es} 

%    General info
\subjclass{Primary 11Y50, 65A05; Secondary 11D25, 14H52}

\keywords{Hall's conjecture, Mordell's equation}

\begin{abstract}
In relation to \Hc, a new algorithm is presented to search for small nonzero $k=\abs{x^3-y^2}$ values. Seventeen new values of $k<x^{1/2}$ are reported.
\end{abstract}

\maketitle

\section{Hall's conjecture}
Dealing with natural numbers, the difference 
\begin{equation}\label{Mordell}
k=x^3-y^2
\end{equation}
 is zero when $x=t^2$ and $y= t^3$ but, in other cases, it seems difficult to achieve small absolute values. For a given $k\ne 0$, (\ref{Mordell}) known as Mordell's equation, is an elliptic curve and has only finitely many solutions in integers by Siegel's theorem. Therefore, for any $k \ne 0$, there are only finitely many solutions in $x$ (which is hence bounded). There is a proved bound due to A. Baker \cite{Baker} and improved by H. M. Stark \cite{Stark}, that places the size of $k$ above the order of $\log^c(x)$ for any $c < 1$. 
A bound concerning the minimal growth rate of $|k|$ was found early by M. Hall \cite{BCHS,Hall} by means of a parametric family of the form
\begin{equation}\label{polynomial}
\begin{array}{rcl}
\vspace{1mm}
f(t)  &=& \displaystyle \frac t 9 (t^9+6t^6+15t^3+12), \\ 
\vspace{1mm}
g(t)  &=& \displaystyle \frac {t^{15}} {27} + \frac{t^{12}+4t^9+8t^6} 3 + \frac {5t^3+1} 2, \\  
f^3(t) - g^2(t) &=& - \displaystyle \frac {3t^6+14t^3+27} {108},
\end{array}
\end{equation}
with $t$ congruent to 3 mod 6, which supplies infinitely many cases with $k < Cx^{3/5}$, where $C$ is a positive constant.

H. Davenport \cite{Davenport}, pointed out that the degree of $f^3(t) - g^2(t)$  is always greater than the half of the degree of $f(t)$. This fact and experimental results for $x < 700.000$, prompted Marshall Hall \cite{Hall} to conjecture that $\abs k$ can not be less than $Cx^{1/2}$, for some constant $C$ whose tentative value was fixed to be $C=1/5$. 
Later, L. V. Danilov \cite{Danilov} found a infinite family derived from the unbounded solutions of quadratic equations supplying values of $\abs k < 217 \sqrt 2 x^{1/2}$. N. D. Elkies \cite{Elkies} revised and improved the method reporting the Fermat-Pell family
\begin{equation}\label{Fermat-Pell}
(5^5t^2+3000t+719)^3-(5^3t^2-114t+26)(5^6t^2-5^3123t+3781)^2 = 27(2t-1),
\end{equation}
for $t$ such that $5^3t^2-114t+26$ is a perfect square. This family yields infinitely many solutions with $\abs k < 0.966\,x^{1/2}$.
Theoretical probabilistic considerations (see S. Lang \cite[p. 163]{Lang}) and the finding by computer search of smaller values of $\abs k$ as the present record found by Elkies \cite{Elkies} of $\abs k = 1/46.6\, x^{1/2}$ suggest that 1/2 is not the better exponent and the original conjecture was reformulated in a weaker way as follows: {\it For any exponent $e < 1/2$, a constant $C_e > 0$ exists such that $|x^3-y^2| > C_e x^e$}. \Hc is considered a particular case of the related and more general ABC conjecture \cite{Oesterle,Lang} and both seem  hard to be proved or disproved. S. Mohit and M. R. Murty \cite{Mohit} show that Hall's conjecture implies that there are infinite many primes such that $a^{p-1} \not\equiv 1 \pmod {p^{16}}$ for any $a$.

Since, at present, \Hc is neither proved nor disproved, it is worthy to enumerate the cases when $|k| < x^{1/2}$, what we address as {\it good examples of Hall's conjecture} borrowing notation used in ABC conjecture. Table 1, gathers the currently known examples (excluding the infinite Danilov-Elkies family) displaying the value of $x$ and $r=\sqrt x/k > 1$. The values of $y$ and $k$ are not displayed because $y$ can be calculated as the nearest integer to $x^{3/2}$. Items \#2, \#3, \#12 and \#13 correspond to (\ref{polynomial}) for $t=-3,3,-9$ and 9 respectively. Items \#4 and \#24 and \#43 correspond to the two first solutions of (\ref{Fermat-Pell}) which, in fact, yields infinitely many examples with $r > 1.035$.  J. Gebel, A. Peth\"o and H.G. Zimmer \cite{GPZ} applied their own algorithm to compute the integral points on elliptic curves in all the cases of (\ref{Mordell}) with $|k| \le 100,000$. Their algorithm was able to find up to the item \#14 in table 1 and excludes the possibility of the existence of other cases with $\abs k \le 100,000$. N.D. Elkies \cite{Elkies} gives an extensive account on the subject and develops a lattice base reduction algorithm finding the next items up to \#26. 

In this paper we present a new algorithm that adds 17 new items to the table of the known {\it good examples of Hall's conjecture}. The algorithm has also found 704 cases with $\sqrt x/k \in (0,16]$ which turn out to be uniformly distributed in this range. The data also show that the number of cases for $x< X$ and $k \le n\sqrt x$ may be estimated as $0.80\,n\log(X)$. The algorithm  was developed on an euristical approach, noticing the patterns (see Fig. 1) of the plots of  solutions of (\ref{Mordell}) such that $\abs k$ is minimal for each $x$ value rather than using basis reduction algorithms, common in this kind of problems. The algorithm searches a family of polynomials that contain those minimal solutions. It has the drawback that it does not find all the cases in increasing order and we can not prove that the algorithm does not miss any solution. Nevertheless, it was able to find all the {\it good examples of \Hc} known before this work, including the first cases of the infinite Danilov-Elkies family (we conjeture that the new algorithm could also find the rest of cases, with enought time of computation).

\section{Values of $k$ for $x$ near the square of a rational}
We consider the smallest $k$ satisfying (\ref{Mordell}) as a function of $x$. Using the notation $\round{x}=\left\lfloor x+ 1/2 \right\rfloor$, it is given by
\begin{equation}\label{k(x)}
k(x) = x^3 - \round{x^{3/2}}^2,
\end{equation}
where $y=\round{x^{3/2}}$. It is clear that $k(x)=0$ for $x=t^2$, $t$ integer. 

These expressions of $x$ and $y$  are somehow related with the
expressions given by N.D.Elkies in equations (38) and (40) of \cite{Elkies}, where
$x$ and $y$ are approximated by multiples of a square and a cube,
respectively. Specifically, he writes $x= 3\alpha ^2 + \beta$ and
$y=6\alpha ^3 +3\alpha\beta +\delta$. Then,  Elkies argues that  the
problem can be studied via lattice reduction. We use a different approach, maybe more heuristic and intuitive, which does not employ lattice reductions.

 Our method starts by allowing rational values for $t$ and considering values of $x$ for the integers close to $t^2$. We define
\begin{equation}\label{t}
t=a/b,\quad \alpha \equiv a^2 \pmod {b^2} \quad \mbox{with}\quad|\alpha| \le b^2/2 \quad \mbox{and}\quad x_0 = \lfloor t^2 \rceil = \frac {a^2 - \alpha} {b^2}.
\end{equation}
Function (\ref{k(x)}) for the integer values of $x$ near $x_0$, say for $x_0+i$, where $i=0, \pm 1, \pm 2, \cdots$, takes the form
\begin{equation}\label{k(x_0+i)}
k(x_0+i)= (x_0+i)^3-y^2,\quad y=\round{(x_0+i)^{3/2}}.
\end{equation}

\begin{lemma}\label{lema}
Function (\ref{k(x)}) for integer values of the variable near $t^2$, $t$ rational, is equivalent to the equation,
\begin{equation}\label{k(C,i)}
k(x_0+i)= \frac {-2Ca^3+\frac 3 4(b^2i-\alpha)^2 a^2 - 3 (b^2i-\alpha)Ca - C^2 +(b^2i-\alpha)^3} {b^6},
\end{equation}
for appropriate values of $C$ and $i$ satisfying
\begin{equation}\label{C}
2a^3 + 3{(b^2i-\alpha)a} + 2C \equiv 0 \pmod {2b^3}
\end{equation}
such that $\abs {k(x_0+i)} \le (x_0+i)^{3/2} + \frac 1 4$.
\end{lemma}
\begin{proof}
We note that, by the definition of $x_0$,
$$
(x_0+i)^3= \frac {a^6 +3(b^2i-\alpha)a^4+3(b^2i-\alpha)^2a^2+(b^2i-\alpha)^3} {b^6}.
$$
Let
$$
y^2 = \frac {a^6 +3(b^2i-\alpha)a^4 + p(a)} {b^6}=\left( \frac {a^3+Ba+C} {b^3}\right)^2
$$
for an appropriate third degree polynomial $p(a)$ and coefficients $B$ and $C$. Since $y$ must be an integer, the following two conditions must be fulfilled,
\begin{enumerate}
\item[(i)] $(a^3+Ba+C)^2=a^6+3(b^2i-\alpha) a^4 + p(a)$. \\
\item[(ii)] $ b^3 \mid a^3+Ba+C.$
\end{enumerate}
Expanding the square of the trinomial in (i) we can deduce that $B=3/2(b^2i-\alpha)$ and $p(a)$ takes the following expression on $C$,
$$
p(a)=2Ca^3+\frac 9 4 (b^2i-\alpha)^2a^2+3(b^2i-\alpha)Ca+C^2.
$$
The condition (ii) becomes now $a^3 + \frac 3 2(b^2i-\alpha) a + C \equiv 0 \pmod {b^3}$. Multiplying by 2, we have (\ref{C}).
Substituting the value of $(x_0+i)^3$ and $y^2$ in (\ref{k(x_0+i)}), we have (\ref{k(C,i)}) where values $C$ and $i$ must fulfill (\ref{C}). The condition  $\abs {k(x_0+i)} \le (x_0+i)^{3/2} + 1/4$ easily follows from the bounds of (\ref{k(x)}).Note that $2C$ is only determined modulo $2b^3$ by (\ref{C}), so that we can take $\abs{2C} < b^3$ and the absolute value of the leading term of (\ref{k(C,i)}) is always smaller than $a^3/b^3 \approx x_0^{3/2}$, hence the bounds of (\ref{k(C,i)}) are the bounds of (\ref{k(x_0+i)}) for moderate $i$ values. 
\end{proof}

The allowed values of $C$ imposed by (\ref{C}) form a small set for a fixed $t$, due to the coefficient $b^2$ of $i$.
Fig. 1 depicts an example for $t=a/b=222272/15$. In this case five values are possible, all of them congruent mod ${3b^2}$ and for each one of the five $C$ values, the points split as integer points of five cubic polynomials. The number of polynomials and the values of $C$ depend from $b$, the parity of $a$ and whether 3 divides $b$ or not. The analysis is done in the next section.
\begin{figure}
\centerline{\psfig{figure=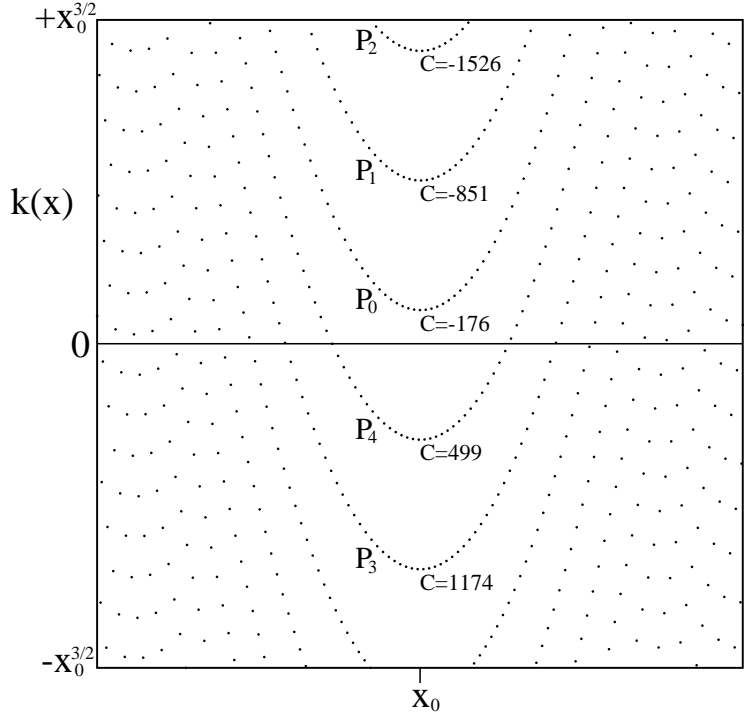,height=10cm}}
\caption{The polynomials for $t=222272/15$}
\end{figure}

\section{The geometry of $k(x)$}
From (\ref{C}), we obtain
$$
3b^2ai \equiv 3 a \alpha - 2a^3 - 2C \pmod {2b^3}.
$$
Both sides of the equation must be divisible by $g=\gcd(3b^2a,2b^3)$ to be solvable in $i$. Note that, since the fraction $t$ is in its lowest terms, $\gcd(a,b)=1$, then 
$$
g = c_1c_2b^2,\;
\left\{
\begin{tabular}{l}
$c_1 = 2$ if $a$ is even, $c_1=1$  otherwise; \\
$c_2 = 3$ if $3 | b$, $c_2=1$  otherwise. \\
\end{tabular} \right. 
$$
Canceling by $g$, we obtain
\begin{equation}\label{ies}
\left( \frac {3a} {c_1c_2} \right) i \equiv\left( \frac {3a\alpha- 2a^3 -2C} {c_1c_2b^2} \right) 
\pmod {\frac {2b} {c_1c_2}}.
\end{equation}
Since the right side of the equation must be an integer, we have that 
$$
3a\alpha- 2a^3 -2C \equiv 0 \pmod {c_1c_2b^2}.
$$
Note that $\alpha = a^2+rb^2$ for some integer $r$, then $3a\alpha \equiv 3a^3+ 3arb^2 \equiv 3a^3 \pmod {c_1c_2b^2}$ for any of the possible values of $c_1$ and $c_2$. Then,
\begin{equation}\label{C_0}
(2C) \equiv a^3 \pmod {c_1c_2b^2}\;
\left\{
\begin{tabular}{l}
$(2C) \equiv a^3  \pmod {c_2b^2}$, when $a$ is odd,\\
$C \equiv a^3/2 \pmod {c_2b^2}$, when $a$ is even,
\end{tabular} \right. 
\end{equation}
The values of $C$ are of the form $C= C_0 + sc_2b^2$, for some integer $s$. Note that whenever $a$ is odd, $C$ may be a rational of the form $C=C'/2,\quad C'$ integer. For $a$ even, $C$ is always integral.

For each value of $C$, the allowed values of $i$ are of the form $i = j + \omega b', \quad b'= 2b/(c_1c_2)$, where $\omega\in\Z$, being $j$ a solution of (\ref{ies}) for such $C$ value.
We have splited up the points $x_0+i$ into $b'$ classes of congruence. For each one, we can associate an unique value of $C$. For each one of the $b'$ couples $(j,C)$, we can substitute the value of $i$ in (\ref{k(C,i)}), getting a cubic polynomial in the variable $\omega$ that contains the family of points corresponding to a value $j$. The effective value of $C$ (out of the $b'$ significant possible values) must be chosen so that $k(x)$ is bounded as imposed by Lemma \ref{lema}. 

{\bf An example:}
We apply the theory to the case illustrated in fig.~1 where $a=222272$ and $b=15$. Then by (\ref{t}),
$$
\alpha=a^2 \mbox{ mod } b^2 = 109,\quad  x_0=(a^2-\alpha)/b^2= 219577075.
$$
Since $a$ is even $c_1=2$ while $c_2=3$ because $3|b$. Applying (\ref{C_0}) we have that
$$
C_0 = a^3/2 \mbox{ mod } 3b^2 = 499, \mbox{ then } C= C_0 + 3b^2s = 499 + 675s, \quad s \in \Z.
$$
Applying (\ref{ies}) we guess that for $C=499$, $i \equiv 4\quad \mbox{mod } 5$, then $i = 4 + 5\omega$. That is, for $x = x_0 + 4 + 5\omega$, the values of $k(x)$ correspond to the points labeled as $P_4$ in the figure.  

\begin{theorem}
The integer points of $k(x) = x^3 - \round{x^{3/2}}^2$, for $x$ near the square of a rational $t=a/b$, are integer points of a set of polynomials
$$
P_j(x)=x^3-(\frac 3 2 t)^2x^2+(\frac 3 2 t^4-  \frac {3C} {b^3}t)x+ \frac C {b^3}t^3- \frac {C^2} {b^6} - \frac 1 4 t^6.
$$
for $0 \le j < b'$, $b'=2b/\gcd(3a,2b)$ and a computable number $C$, provided $|P_j(x)| \le x^{3/2} + 1/4$ and $x=\round{t^2}~+~j+~\omega~b'$, $\omega\in\Z$.
\end{theorem}
\begin{proof}
We note that, for each $C$ value in (\ref{ies}), $i$ is in $Z_{b'}$, with $b'=2b/(c_1c_2)= 2b/\gcd(3a,2b)$. Then, we put $i= j +\omega b'$ and we have $b'$ cases of values $(C,j)$ for each $t$ value. We substitute $j+\omega b'$ for  $i$ in (\ref{k(C,i)}). Then the term $(b^2i -\alpha)$ becomes $(b^2(j+\omega b')-\alpha)=b^2(b'\omega +j- \frac \alpha {b^2})$. We get a polynomial in $\omega$ for each $j$,
$$
P_j(\omega)=(b'\omega+j-\frac \alpha {b^2})^3+\frac 3 4 t^2(b'\omega+j-\frac \alpha {b^2})^2- \frac {3C}{b^3}t(b'\omega+j-\frac \alpha {b^2})- \frac {2C}{b^3}t^3-\frac {C^2}{b^6}.
$$
Recall than $x=x_0+j+b'\omega$. Let be $z=x-t^2$, then
\begin{equation}\label{z}
z=x_0+j+b'\omega - t^2=j+b'\omega - \frac \alpha {b^2}.
\end{equation}
as follows from (\ref{t}). Then
\begin{equation}\label{P_j(z)}
P_j(z)= z^3+\frac 3 4 t^2z^2-\frac {3C} {b^3}tz-\frac {2C} {b^3}t^3 - \frac {C^2} {b^6}.
\end{equation} 
We can put the polynomial in function of $x$ by means of the variable change $z=x-t^2$, obtaining  the polynomials $P_j(x)$ in the statement of the theorem.

As in Lemma 2.1, we must impose that $k(x)$ is bounded by $x^{3/2}+1/4$. We can compute $C$ by (\ref{C_0}).
\end{proof}

We analyze the polynomial for $x$ near $t^2$, i.e. for $z \approx 0$. Note that, for small $z$, the cubic term is small compared with the quadratic one, then, it is locally approximated by a parabola with the vertex near $z=0$. The value of the polynomial at the vertex is therefore near 
$$
P_j(z=0)\approx -\frac {2C}{b^3}t^3.
$$ 
Since the values of $C$ for some $t$ are equally spaced (they are of the form $C=C_0+sc_2b^2$ with $s$ integral), the vertices are also approximately equally spaced. Note that, for $C>0$, the value of the polynomial is negative in the vertex, but the polynomial grows until crossing $P_j(z)=0$. Recall that for integral $\omega$, the polynomials give $k(x)$. In this way, it is possible to find a small $k(x)$, selecting for each fixed $a/b$ the polynomial with lower positive $C$ value (the polynomial $P_4$ in the example of fig.~1) and calculating the value $\omega$ where the polynomial has its smallest absolute value. Nevertheless, better results can be achieved calculating the optimal value of $a$ for each value $b$, for which we can select an optimal polynomial. First, we can impose $C$ to be positive and small (say $C=1/2,1,3/2 \cdots$). We can also impose that the smallest value of the polynomial is for $\omega=0$, so the value of $z$ in (\ref{z}) is small and we can hope a small value of $P_j(z)$. Even more, if $j=0$, then $z=\alpha/b^2 \le 1/2$ which is the smallest value possible for $z$. We can approximate the value of $P_0$, for $\omega=0$ by taking only the terms in $t^3$ and $t^2$ in (\ref{P_j(z)}) getting
\begin{equation}\label{P_0}
P_0(\omega=0) \approx \frac 3 4 t^2 \frac {\alpha^2} {b^4} - \frac {2C} {b^3}t^3. 
\end{equation}
Suppose that $a<O(b^{5/2})$, $\alpha < O(b^{5/4})$ and $C=O(1)$. Then $x\approx a^2/b^2 < O(b^3)$ and both terms of the approximation are smaller than $O(x^{1/2})$. Even if $a$ and $\alpha$ are above those bounds, the polynomial may have a small value if both terms are balanced.

\section{The algorithm}
The algorithm runs over $b=2,3,4,\dots,$ and for small values of $C$ whose bounds will be discussed later. For each pair $(b,C)$, we must compute the values of $a$ (if any value exists) in (\ref{C}) corresponding to the polynomial $P_0$ (that is for $j=0$) whose smallest absolute value is for $\omega=0$. So, (\ref{C}) becomes
\begin{equation}\label{Cj0}
2a^3 - 3\alpha + 2C \equiv 0 \pmod {2b^3}.
\end{equation}
First, we compute the value of $a$ modulo $b^2$ and then we lift to modulo $2b^3$. From (\ref{C_0}), we have that 
\begin{equation}\label{a_0}
a \equiv (2C)^{1/3}  \pmod {b^2}.
\end{equation}
Note that the equation is solvable only when $\gcd(2C,b)=1$ because $a$ and $b$ are coprime. The cube root may be computed, for example, using the algorithms described in \cite{German}.
Before lifting the solution of (\ref{a_0}) to modulo $2b^3$ in (\ref{C}), we must analyze the possible values of the variables whose constraints can be summarized as
\vskip 5mm
$
\gcd(2C,b)=1,\quad b
\left\{
\begin{tabular}{l} 
  even, $C$ rational, \\
  odd, $\left\{ \begin{tabular}{l} 
        $C$ rational; \\
        $C$ integer.
        \end{tabular} \right.$
\end{tabular} \right.
$
\vskip 5mm
In fact, if $b$ is even, $a$ must be odd because $t$ is in its lower terms.  Being $b$ even and $a$ odd, $\alpha$ must be odd, by the definition of $\alpha$ in (\ref{t}), what implies that $(2C)$ in (\ref{Cj0}) must be odd. Therefore, $C=C'/2$ for some odd integer $C'$. If $b$ is odd, $a$ can be even or odd, then $C$ may be rational or integer.

Let $a_0$ be a solution of (\ref{a_0}), then $a=a_0+kb^2$ for some integer $k$ that must fulfill (\ref{Cj0}). So, we have
$$
2(a_0+kb^2)^3 - 3\alpha (a_0+kb^2) + 2C \equiv 0 \pmod {2b^3}.
$$
Expanding and neglecting the vanishing terms, we have
$$
2a_0^3 + 6a_0^2kb^2 - 3\alpha a_0 - 3\alpha kb^2 + 2C \equiv 0 \pmod {2b^3}.
$$
Then
$$
3b^2(2a_0^2 - \alpha )k \equiv - \left( 2a_0^3 -  3\alpha a_0 + 2C \right) \pmod {2b^3}.
$$
Let 
$$ 
d = \gcd(3(2a_0^2 - \alpha ),2b).
$$
Then
\begin{equation}\label{k}
k \equiv - \left( \frac {3(2a_0^2 - \alpha )} d \right)^{-1}  \left( \frac {2a_0^3 -  3\alpha a_0 + 2C}  {db^2} \right) \pmod { \frac {2b} d}.
\end{equation} 
 Note that the possible values of $d$ may be $d=1$, $d=3$ if $ 3|b$, $d=2$ if $\alpha$ is even or $d=6$ if both conditions holds.
Let $k_0$ be a solution of the above equation, then $k=k_0+n2b/d$ for some integer $n$. Substituting the value of $k$ in $a=a_0+kb^2$, we have finally
$$
a = a_0 +k_0b^2 + n \frac {2b^3} d.
$$
Since we are looking for small values of $P_0$ in (\ref{P_0}), that is for $\omega=0$, we want that $a \approx (3\alpha^2)/(8C)$. Then
$$
n=\round{\frac d {2b^3}\left(\frac {3\alpha^2} {8C} -a_0-k_0b^2\right)}.
$$
Once the most appropriate values $a$ for each couple $(b,C)$ are computed, we can compute $x_0=(a^2-\alpha)/b^2$  expecting to find a small value of $k(x)$.

\section{Results and conclusions}
The algorithm was programed as a PARI script and translated to C with the utility GP2C \cite{PARI2}. The executable was run in various computers during several periods being the total time of computation equivalent to about 1000 days in a AMD-64 3400+. The algorithm was able to find all the so far known items and the remaining 17 new items displayed in table 1 with $r = \sqrt{x}/k \ge 1$. Three items were found by Johan Bosman using the software of the authors. 

\begin{center}
\begin{tabular}{|l|l|r|l|r|l|} \hline
\# & \multicolumn{1}{c|} {$x$}  &\multicolumn{1}{c|} {$r$} & \multicolumn{1}{c|}{$b$} & \multicolumn{1}{c|}{$C$}  & Comments \\  \hline \hline
1 & 2 & 1.41 & - & - & -\\
2 & 5234 & 4.26 & 26 & 1/2  &H, GPZ, (\ref{polynomial})\\
3 & 8158 & 3.76 & 28 & 1/2  &H, GPZ, (\ref{polynomial})\\
4 & 93844 &1.03 & 53 & 1   &H, GPZ, (\ref{Fermat-Pell})\\
5 & 367806 & 2.93 & 117 & 1/2   &H, GPZ \\
6 & 421351 & 1.05 & 26 & 1/2  & H, GPZ\\
7 & 720114 & 3.77 & 42 & 1/2  & H, GPZ \\
8 & 939787 & 3.16 & 115 & 2  & H, GPZ\\
9 & 28187351 & 4.87 & 159 & 5  & H, GPZ\\
10 & 110781386 &  1.23 & 95 & 1/2  & H, GPZ\\
11 & 154319269 &  1.08 & 228 & 1/2  & H, GPZ\\
12 & 384242766 & 1.34 & 728 & 1/2  & H, GPZ, (\ref{polynomial})\\
13 & 390620082 &  1.33 & 730 & 1/2  & H, GPZ, (\ref{polynomial})\\
14 & 3790689201 & 2.20 & 1155 & 4  & GPZ\\
15 & 65589428378 & 2.19 & 5235 & 17/2  & E \\
16 & 952764389446 & 1.15 & 1448 & 5/2  &  E \\
17 & 12438517260105 & 1.27 & 13415 & 6  & E \\
18 & 35495694227489 & 1.15 & 97266 & 1/2  & E \\
19 & 53197086958290 & 1.66 & 13777 & 1  & E \\
20 & 5853886516781223 & 46.60 & 137035 & 9  & E, !\\
21 & 12813608766102806 &  1.30 &  6291 & 35/2  & E \\
22 & 23415546067124892 & 1.46 & 1315447 & 32  & E, * \\
23 & 38115991067861271 & 6.50 & 321346 & 1/2  & E \\
24 & 322001299796379844 & 1.04 & 1313479 & 11 & E, (\ref{Fermat-Pell}) \\
25 & 471477085999389882 & 1.38 & 3281374 & 95/2 & E \\
26 & 810574762403977064 & 4.66 & 5346121 & 49/2  & E \\
27 & 9870884617163518770 & 1.90 & 4928788 & 109/2  & JHS \\
28 & 42532374580189966073 & 3.47 & 583876 & 9/2  & JHS\\
29 & 51698891432429706382 & 1.75 &     19061951 & 29  & JHS \\
30 & 44648329463517920535 & 1.79 &     11744301 & 13   & JHS \\
31 & 231411667627225650649 & 3.71 &    11694866 & 347/2  & JHS \\
32 & 601724682280310364065 & 1.88 &    7496613 & 13  & JHS \\
33 & 4996798823245299750533 & 2.17 &   76010518 & 67/2  & JHS \\
34 & 5592930378182848874404 &  1.38  & 93203798 & 139/2 & JHS \\
35 & 14038790674256691230847 & 1.27 &  61769318 & 53/2 & JHS \\
36 & 77148032713960680268604 & 10.18 & 184388019 & 4 &  JB \\
37 & 180179004295105849668818 & 5.65 & 292889921 & 45/2 &  JB \\
38 & 372193377967238474960883 & 1.33 & 2554989 & 4 & JHS \\
39 & 664947779818324205678136 &  16.53  & 678534061 & 39.2 & JHS \\
40 & 2028871373185892500636155 & 1.14 &  490670918 & 55/2 & JB \\
41 & 37223900078734215181946587 & 1.87  & 530793746  & 1457/2  & JHS \\ 
42 & 3690445383173227306376634720 & 1.51  & 685266726  & 1/2  & JHS \\
43 & 1114592308630995805123571151844 & 1.04  & 52019836686   & 1475/2 & (\ref{Fermat-Pell}) \\ 
44 & 6078673043126084065007902175846955 & 1.03  & 8144029787  & 3 & JHS \\

\hline
\end{tabular}
\end{center}
%\vskip -.2cm
\noindent{\bf Table 1: Good examples of \Hc}.\newline
{\bf H}; Found by M. Hall \cite{Hall}.
{\bf GPZ}: Found by J. Gebel, A. Peth\"o and H. G. Zimmer \cite{GPZ} with the elliptic curve method. \newline
{\bf E}: Found by N.D. Elkies with a lattice basis reduction algorithm.\newline
{\bf JHS}: Found by the authors of this paper.\newline
{\bf JB}: Found by Johan Bosman with the software of the authors.\newline
{\bf (\ref{polynomial})}: Items \#2,\#3,\#12 and \#13 correspond to the polynomial family for $t=-3,3,-9$ and 9 respectively in equation (\ref{polynomial}).\newline
{\bf (\ref{Fermat-Pell})}: The two first terms of the infinite Fermat-Pell family from equation (\ref{Fermat-Pell}).\newline
{\bf !}: Record of N.D. Elkies. \newline
{\bf *}: Obtained from \#20 scaling $(x,y,k)$ to $(2^2x,2^3y,2^6k)$.
\vskip .4cm

The value of the integer $b$ runs from 2 to a bound $B$ while for each value of $b$, values of $C$ from 1/2 to $b^u$ were checked. Experimentally, it was found that a good interval for $C$ correspond to $u=1/3$. With a shorter interval, for example $u=1/4$, the algorithm runs faster but can not find some items in the table. In the first case, for $1/2 \le C \le b^{1/3}$ the bound for $b$ reached was $B=6\cdot10^8$, while for $1/2 \le C \le b^{1/4}$ the bound was $B=5\cdot10^9$. Item \#42 and \#43 was found constraining the search to lower values of $C$. Although the used bounds for $C$ seem reasonable, it can not be proved that the algorithm can find all existing {\it good examples of Hall's conjecture}.

The exact complexity of the algorithm must be expressed only in function of the bound of $b$. For each $b$ value, we must factorize $b$ and, for each of the $b^u$ values of $C$, we must perform a cube root and other polynomial time operations. Suppose that the factorization of $b$ takes $O(b^v)$ operations. Then it takes $O(B^{ 1+\max \left\{u,v \right\}}\log^{O(1)}(B))$ operations to explore the $b\le B$ values. Since  a simple trial division with a table of primes was used, the exponent $v$ must be under 0.35 in about half of the instances (see \cite{Knuth}). It is not possible to express the complexity as a fixed function of the bound of $x$, nevertheless a estimation can be performed. The value of $a$ is $O(b^\theta)$ for an variable exponent $\theta$ that is 2.5 in average (see the last paragraph in Section 3). Then, the corresponding value of $x$ is $x\approx a^2/b^2 = b^{2\theta-2}$. Since $\max\left\{u,v\right\}=0.35$, the complexity takes the form $O(X^{\frac {1.35} {2\theta-2}} \log^{O(1)}(X))$ for $x < X$. If we consider the average value $\theta=2.5$ we obtain an equivalent complexity of $O(X^{0.45})$ which coincides with the experimental average complexity.

Another search was performed seeking for cases in which $k \le 16\sqrt{x}$ for $b < 12,500,000$. The algorithm discovered 704 cases whose analysis shows that the distribution of the values of $k/\sqrt x \in (0,16]$ is very close to a uniform distribution for this range of $k$ (the mean is 7.996 and the Kolmogorov-Smirnov test gives a confidence of 98.60\% to the hypothesis of a uniform distribution). This gives extra support to the possibility that the size of $k/\sqrt x$ may be arbitrarily small, and so the corresponding exponent in Hall's conjecture must be at most $1/2-\epsilon$. The data also support that the number of cases for $x< X$ and $k \le n\sqrt x$ may be estimated as $0.80\,n\log(X)$, somewhat less than the expected value $n\log(X)$.

\bibliographystyle{amsplain}

\end{document}